\documentclass[review]{elsarticle}

\usepackage{tikz}
\usetikzlibrary{decorations.pathmorphing}
\tikzset{every node/.style={circle,fill, inner sep=0pt, minimum size=1.5mm}, every picture/.style={scale=.85}}

\usepackage{hyperref}

\usepackage{amsmath,amssymb,amsthm}
\allowdisplaybreaks

\newcounter{tbox}

\newcommand{\sta}[1]{\vspace*{0.3cm}\refstepcounter{tbox}\noindent{ \parbox{\textwidth}{(\thetbox) \emph{#1}}}\vspace*{0.3cm}}

\newtheorem{theorem}{Theorem}
\newtheorem{lemma}[theorem]{Lemma}

\journal{Discrete Mathematics}

  \newcommand{\joe}[1]{{#1}}

\begin{document}
	
	\begin{frontmatter}

\joe{
	\title{Minimal induced subgraphs of two classes of 2-connected non-Hamiltonian graphs}
}
		
		\author{Joseph Cheriyan\fnref{jc}}
		
		\author{Sepehr Hajebi} 
		
		\author{Zishen Qu} 
		
		\author{Sophie Spirkl\fnref{sophie}}
		
		\fntext[jc,sophie]{ We acknowledge the support of the Natural Sciences and Engineering Research Council of Canada (NSERC), [funding reference numbers RGPIN-2019-04197 and RGPIN-2020-03912]. Cette recherche a été financée par le Conseil de recherches en sciences naturelles et en génie du Canada (CRSNG), [numéros de référence RGPIN-2019-04197 et RGPIN-2020-03912].}
		\address{Department of Combinatorics and Optimization, University of Waterloo, Waterloo, Ontario, N2L 3G1, Canada}

\begin{abstract}
In 1981, Duffus, Gould, and Jacobson showed that every connected graph either has a Hamiltonian path, or contains a claw ($K_{1,3}$) or a net (a fixed six-vertex graph) as an induced subgraph. This implies that subject to being connected, these two are the only minimal (under taking induced subgraphs) graphs with no Hamiltonian path. 

Brousek (1998) characterized the minimal graphs that are
$2$-connected, non-Hamiltonian and do not contain the claw as an induced subgraph. We characterize the minimal graphs that are $2$-connected and non-Hamiltonian for two classes of graphs: (1) split graphs, (2) triangle-free graphs. We remark that testing for Hamiltonicity is \textsf{NP}-hard in both classes.
\end{abstract}

	\end{frontmatter}
		

	\section{Introduction}
	
	Graphs in this paper are finite and without loops or parallel edges. For a graph $G$ and $X \subseteq V(G)$, $G[X]$ denotes the induced subgraph of $G$ with vertex set $X$, and $G \setminus X$ denotes $G[V(G) \setminus X]$. A \emph{Hamiltonian path} (resp.\ \emph{Hamiltonian cycle}) in a graph $G$ is a (not necessarily induced) subgraph $H$ of $G$ which is a path (resp.\ cycle), and $V(H) = V(G)$. A graph is \emph{Hamiltonian} if it has a Hamiltonian cycle. 
	
	We say that a graph $H$ is an \emph{HP-obstruction} if $H$ is connected, has no Hamiltonian path, and every induced subgraph of $H$  either equals $H$, or is not connected, or has a Hamiltonian path. Analogously, a graph $H$ is an \emph{HC-obstruction} if $H$ is $2$-connected, has no Hamiltonian cycle, and every induced subgraph of $H$  either equals $H$, or is not $2$-connected, or has a Hamiltonian cycle. 
	
	The \emph{claw} is the complete bipartite graph $K_{1,3}$. The \emph{net} is the unique graph with degree sequence $(3,3,3,1,1,1)$, and equivalently the graph with vertex set $\{a,b,c, a', b', c'\}$ and edge set $\{a'b', b'c', a'c', aa', bb', cc'\}$. The \emph{snare} is the graph obtained from a net by adding a vertex and making it adjacent to every vertex of the net. The following theorem of Duffus, Gould, and Jacobson characterizes all HP-obstructions.
	
	\begin{theorem}[\cite{duffus-gould-jacobson}; see also \cite{shepherd}]
	    There are exactly two HP-obstructions: the claw and the net. 
	\end{theorem}

	Following the same line of thought, in this note we are interested in understanding HC-obstructions. 	
\joe{
In \cite{brousek}, Brousek gave a complete characterization of HC-obstructions that do not contain the claw as an induced subgraph, and Chiba \& Furuya \cite{chiba-furuya} further studied induced subgraphs of non-minimal 2-connected non-Hamiltonian graphs. Ding \& Marshall \cite{ding-marshall} obtained a complete characterization in the case when ``induced subgraph'' is replaced by ``induced minor'' in the definition of an HC-obstruction.  
}

Let us describe our main results. A \emph{clique} in a graph $G$ is a set $K$ of pairwise adjacent vertices. A \emph{stable set} in a graph $G$ is a set $S$ of pairwise non-adjacent vertices. A \emph{split graph} is a graph $G$ with a partition $(S, K)$ of $V(G)$ such that $S$ is a stable set and $K$ is a clique in $G$.
		
		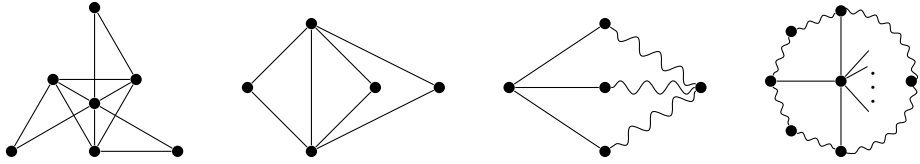
\begin{figure}[h!]\centering
			\begin{tikzpicture}
			\newdimen\R
			\R=0.75cm
			\newdimen\br
			\br=1.5cm
			\node (c2) at (270:\R) {};
			\node (b2) at (150:\R) {};
			\node (a2) at (30:\R) {};
			\node (p) at (0,0) {};
			\node (a1) at (90:\br) {};
			\node (b1) at (210:\br) {};
			\node (c1) at (330:\br) {};
			\draw (p) -- (c2);
			\draw (p) -- (a2);
			\draw (p) -- (b2);
			\draw (p) -- (c1);
			\draw (p) -- (a1);
			\draw (p) -- (b1);
			\draw (a2) -- (b2) -- (c2) -- (a2);
			\draw (a1) -- (a2);
			\draw (b1) -- (b2);
			\draw (c1) -- (c2);
			\end{tikzpicture}
			\hfill
			\begin{tikzpicture}
			\node (s1) at (-1,0) {};
			\node (k1) at (0,1) {};
			\node (k2) at (0,-1) {};
			\node (s2) at (1,0)	{};
			\node (s3) at (2,0) {};
			\draw (k1) -- (k2);
			\draw (s1) -- (k1);	
			\draw (s1) -- (k2);
			\draw (s2) -- (k1);
			\draw (s2) -- (k2);	
			\draw (s3) -- (k1);	
			\draw (s3) -- (k2);		
			\end{tikzpicture}
			\hfill
			\begin{tikzpicture}
			\node (u) at (-1.5,0) {};
			\node (p1) at (0,1) {};
			\node (p2) at (0,0) {};
			\node (p3) at (0,-1) {};
			\node (v) at (1.5,0)	{};
			\draw (u) -- (p1);
			\draw (u) -- (p2);
			\draw (u) -- (p3);
			\draw[decorate, decoration=snake] (v) -- (p1);
			\draw[decorate, decoration=snake] (v) -- (p2);
			\draw[decorate, decoration=snake] (v) -- (p3);
			\end{tikzpicture}
			\hfill
			\begin{tikzpicture}
			\newdimen\R
			\R=1.1cm
			\node (h2) at (90:\R) {};
			\node (h5) at (270:\R) {};
			\node (ell) at (180:\R) {};
			\node (ell2) at (0:\R) {};
			\node (m1) at (135:\R) {};
			\node (m2) at (225:\R) {};
			\node (c) at (0,0) {};
			\draw (c) -- (h2);
			\draw (c) -- (h5);
			\draw (c) -- (ell);
			\draw [decorate, decoration={snake, segment length=2mm, amplitude=0.35mm}] (h2) -- (m1);
			\draw [decorate, decoration={snake, segment length=2mm, amplitude=0.35mm}] (m1) -- (ell);
			\draw [decorate, decoration={snake, segment length=2mm, amplitude=0.35mm}] (ell) -- (m2);
			\draw [decorate, decoration={snake, segment length=2mm, amplitude=0.35mm}] (m2) -- (h5);
			\draw [decorate, decoration={snake, segment length=2mm, amplitude=0.35mm}] (0,\R) arc(90:-90:\R);
			\newdimen\rstp
			\rstp=0.5cm
			\node[style={fill=none}] (r1) at (\rstp,\R/2) {};
			\node[style={fill=none}] (r2) at (\rstp,\R/4) {};
			\node[style={fill=none}] (r3) at (\rstp,-\R/2) {};
			\draw (c) -- (r1);
			\draw (c) -- (r2);
			\draw (c) -- (r3);
			\path (r2) -- (r3) node [style={fill=none}, font=\large, midway, sloped] {$\dots$};
			\end{tikzpicture}\caption{From left to right: the snare, the $2$-nova, a theta, and a triangle-free wheel. Squiggly edges represent paths of length at least one.}
		\end{figure}
	An \emph{$n$-sun} is a graph obtained from a cycle $C$ with $2n$ vertices $v_1, \dots, v_{2n}$ that occur in this order along $C$ by adding all edges $v_{2i}v_{2j}$ for distinct $i, j \in \{1, \dots, n\}$. An \emph{$n$-nova} is obtained from an $n$-sun by adding a vertex $w$ and edges $wv_{2i}$ for all $i \in \{1, \dots, n\}$. Our first theorem, the following, gives a complete characterization of HC-obstructions that are split graphs.
	
	\begin{theorem} \label{thm:mainsplit}
	    The snare and all $n$-novae for $n \geq 2$ are HC-obstructions. Moreover, these are the only HC-obstructions which are split graphs. 
	\end{theorem}
	
	A \emph{theta} is a graph consisting of two non-adjacent vertices $u$ and $v$ and three paths $P_1, P_2, P_3$ from $u$ to $v$ and each of length at least two, such that the sets $V(P_1) \setminus \{u, v\}, V(P_2) \setminus \{u, v\}, V(P_3) \setminus \{u, v\}$ are disjoint and have no edges between them. The vertices $u$ and $v$ are the \emph{ends} of the theta. A \emph{closed theta} is a graph obtained from a theta with ends $u, v$ by adding the edge $uv$. 
	
	A graph is \emph{triangle-free} if it contains no three-vertex clique. A \emph{wheel} is a pair $(W, v)$ such that $W$ is a cycle, and $v$ is a vertex with at least three neighbours in $W$\footnote{In a standard definition of a wheel, the cycle $W$ is required to be of length at least four. Note that this does not matter for our purposes as we are only concerned with triangle-free wheels.}. 
	
	\begin{theorem} \label{thm:maintf}
	    All thetas, triangle-free closed thetas, and triangle-free wheels are HC-obstructions, and they are the only HC-obstructions which are triangle-free. 
	\end{theorem}

	\section{Split graphs}

	In this section we prove Theorem \ref{thm:mainsplit}. The following is well-known (see, for example, \cite{chvatal}): 
	\begin{lemma} \label{lem:tough}
	    Let $G$ be a graph and $X \subseteq V(G)$. If $G \setminus X$ has more than $|X|$ connected components, then $G$ has no Hamiltonian cycle. 
	\end{lemma}
	
	From this, we deduce:
	\begin{lemma} \label{lem:sg1}
	    The snare and all $n$-novae for $n \geq 2$ are $2$-connected graphs with no Hamiltonian cycle. 
	\end{lemma}
	\begin{proof}
\joe{
	Clearly, these graphs are $2$-connected.
First, consider the snare. Suppose it has a Hamiltonian cycle; then,
a Hamiltonian path of a net can be obtained by deleting one particular
vertex of the snare. This is a contradiction.
Next, consider an $n$-nova for $n \geq 2$.  The graph is non-Hamiltonian,
by Lemma~\ref{lem:tough} with $X = \{v_{2i} : i \in \{1, \dots, n\}\}$, where
the vertex labels are as in the definition.
}
	\end{proof}

	In view of Lemma \ref{lem:sg1}, in order to prove Theorem \ref{thm:mainsplit}, it is sufficient to prove: 
	\begin{theorem}
	    Let $G$ be a $2$-connected split graph with no induced subgraph isomorphic to the snare or an $n$-nova for $n\geq 2$. Then $G$ has a Hamiltonian cycle. 
	\end{theorem}
	\begin{proof}
		Suppose for a contradiction that $G$ has no Hamiltonian cycle. Let $(S, K)$ be a partition of $V(G)$ such that $S$ is a stable set, $K$ is a clique, and subject to this, $|K|$ is maximized.  Then $S \neq \emptyset$, because every $2$-connected complete graph has a Hamiltonian cycle. 
\joe{
\big(Our choice of $(S, K)$ ensures that for each vertex $s\in{S}$,
we have $N(s)\subsetneq{K}$; otherwise, if some vertex $s_0\in{S}$
has $N(s_0)=K$, then $(S \setminus \{s_0\}, K \cup \{s_0\})$ is a
partition of $V(G)$ that contradicts our choice of $(S, K)$.\big)
}

We first prove:
		
		\sta{\label{lt2nbrs} There is a $k \in K$ with $|N(k) \cap S| \geq 3$.}
		
		Suppose not; that is, suppose that each $k \in K$ has at most 2 neighbours in $S$. For every $s \in S$, let us pick two distinct edges $e_1(s)$ and $e_2(s)$ incident with $s$, subject to the number of cycles in $H = (V(G), \{e_i(s) : s \in S, i \in \{1,2\}\})$ being as small as possible. 
		
		If $H$ has no cycles, then, since every vertex in $H$ has degree at most two, $H$ is a disjoint union of a set $K^* \subseteq K$ of isolated vertices and of paths $P_1, \dots, P_t$, each containing at least one edge, with ends in $K$, and with $S \subseteq V(P_1) \cup \dots \cup V(P_t)$. 
		Let $P^*$ be a path containing all vertices of $K^*$ (possibly empty). Then, since all non-empty paths among $P_1, \dots, P_t, P^*$ have ends in $K$ (and thus are pairwise adjacent), the concatenation of $P_1, \dots, P_t, P^*$ is a Hamiltonian cycle of $G$, a contradiction. 
		
\joe{
	Therefore, we may assume that $H$ has a cycle. Since $H$ is not a
Hamiltonian cycle of $G$, it follows that there is a cycle $C$ in $H$
such that $V(C) \neq V(G)$. Since every vertex in $K$ has at most
two neighbours in $S$, it follows that $G[V(C)]$ is a $(|V(C)|/2)$-sun
and  $N(s) \cap V(C) = \emptyset$ for all $s \in V(G) \setminus
V(C)$; therefore, there is a vertex $k' \in K \setminus V(C)$. If
$N(k') \cap V(C) \cap S = \emptyset$, then $G[V(C) \cup \{k'\}]$
is a $(|V(C)|/2)$-nova, a contradiction. Now, let $s \in N(k') \cap
V(C) \cap S$, and consider the graph $H' = (V(G), (E(H) \setminus
\{e_1(s)\}) \cup {sk'})$ obtained by choosing $sk'$ as $e_1(s)$ instead.
%
%
Note that $H'[V(C)]$ is a path (obtained from $C$ by removing the
edge $e_1(s)$), and therefore $H'[V(C) \cup \{k'\}]$ is connected,
contains $sk'$, and contains a vertex of degree one (the end of
$e_1(s)$ in $K$). Since all vertices of $H'$ have degree at most
two, it follows that the component of $H'$ containing $sk'$ is a
path. This contradicts our choice of $H$, since $H'$ has fewer cycles
than $H$, and \eqref{lt2nbrs} follows.
}
		
		\sta{\label{2nbrs}If $s, s' \in S$ are distinct and $|N(s) \cap N(s')| \geq 2$, then $N(s) \cup N(s') = K$.}
		
		Let $k, k' \in N(s) \cap N(s')$ be distinct; and suppose for a contradiction that there is a vertex $k'' \in K \setminus (N(s) \cup N(s'))$. Then $G[\{k, k', k'', s, s'\}$] is a 2-nova, a contradiction. This proves \eqref{2nbrs}. 
		
		\sta{\label{3nbrs}Let $k \in K$ and let $s_1, s_2, s_3 \in S \cap N(k)$ be distinct. Then $N(s_1) \cap N(s_2) \cap N(s_3) = \{k\}$ and $N(s_i) \cup N(s_j) = K$ for all distinct $i, j \in \{1, 2, 3\}$.}
		
		If there is a vertex $k' \in (N(s_1) \cap N(s_2) \cap N(s_3)) \setminus \{k\}$, then $G[\{k, k', s_1, s_2, s_3\}]$ is a 2-nova in $G$, a contradiction; this proves the first part of \eqref{3nbrs}. 
		 
\joe{
Since $G$ is $2$-connected, we may choose $k_i \in N(s_i) \setminus
\{k\}$ for all $i \in \{1, 2, 3\}$.  Note that $k_1, k_2, k_3$ need
not be distinct.  Since $G[\{k, k_1, k_2, k_3, s_1, s_2, s_3\}]$
is not a snare, it follows that $G$ contains an edge  $s_i k_j$ for
distinct $i, j \in \{1, 2, 3\}$. By symmetry, we may assume that
$i = 2$ and $j = 1$. By \eqref{2nbrs}, it follows that $N(s_1) \cup
N(s_2) = K$. By symmetry, we may assume that $k_3 \in N(s_1)$, and
so again by \eqref{2nbrs}, it follows that $N(s_1) \cup N(s_3) =
K$. If $(N(s_2) \cap N(s_3)) \setminus \{k\} \neq \emptyset$, then
\eqref{3nbrs} follows from another application of \eqref{2nbrs}.
Now, assume that $N(s_2) \cap N(s_3) = \{k\}$.
Since $K \setminus N(s_j) \subseteq N(s_1)$ for $j = 2,3$, we have
$K \setminus (N(s_2) \cap N(s_3)) = K \setminus \{k\} \subseteq N(s_1)$.
Therefore, $N(s_1) = K$.  This contradicts our choice of $(S, K)$,
and \eqref{3nbrs} follows.
}

		\sta{\label{4nbrs}$|N(k) \cap S| \leq 3$ for all $k \in K$.}

\joe{
Suppose not; let $k \in K$, and let $s_1, s_2, s_3, s_4$ in $N(k)
\cap S$ be distinct. Let $k' \in N(s_4) \setminus \{k\}$.
Note that at least two of $s_1, s_2, s_3$, say $s_1, s_2$, are adjacent to
$k'$, as otherwise there exist distinct $i, j \in
\{1, 2, 3\}$ with $k' \not\in N(s_i)$ and $k' \not\in N(s_j)$, and
so $N(s_i)\cup N(s_j)\not=K$, a contradiction with~\eqref{3nbrs}. But then $\{k,k'\}\subseteq N(s_1)\cap N(s_2)\cap N(s_4)$, which again violates \eqref{3nbrs}.
This proves \eqref{4nbrs}. \medskip
}

\joe{
Let $k, s_1, s_2, s_3$ be as in \eqref{3nbrs}. For $i \in \{1,2,3\}$, let $K_i = K \setminus N(s_i)$, and note that $K_i$
is non-empty since $N(s_i)\not=K$ (by our choice of $(S,K)$).
It follows that $K_i \subseteq N(s_j)$ for all distinct $i, j \in
\{1,2,3\}$, and that $K = \{k\} \cup K_1 \cup K_2 \cup K_3$.
}

		\sta{\label{s4}$|S| = 3$.}
		
		Suppose that there is a vertex $s_4 \in S \setminus \{s_1, s_2, s_3\}$. If $N(s_4) \cap K_i$ contains two distinct vertices $k', k''$ for some $i \in \{1, 2, 3\}$, then $G[\{k', k''\} \cup \{s_j: j \in \{1,2,3,4\} \setminus \{i\}\}]$ is a 2-nova, a contradiction.
		
		Let $k_1, k_2$ be two distinct neighbours of $s_4$; then $k_1, k_2 \neq k$ by \eqref{4nbrs}, so we may assume by symmetry that $k_i \in K_i$ for $i=1,2$. Now, let $k_3 \in K_3$. Since $G[\{k_1, k_2, k_3, s_4, s_3\}]$ is not a 2-nova, it follows that $s_4$ is adjacent to $k_3$, and therefore to every vertex in $K_3$; and thus, by symmetry, to every vertex in $K_1 \cup K_2 \cup K_3$. This implies that $|K_1| = |K_2| = |K_3| = 1$, since we proved that $s_4$ has at most one neighbour in each of these sets. Also, every vertex in $K$ has three neighbours in $\{s_1, s_2, s_3, s_4\}$, and so by \eqref{4nbrs}, $|S| = 4$. Thus $|V(G)| = 8$ and traversing the vertices in the order $k, s_1, k_2, s_3, k_1, s_4, k_3, s_2$ is a Hamiltonian cycle, a contradiction. This proves \eqref{s4}.\medskip
		
		Now, let $k_i \in K_i$ for $i \in \{1,2,3\}$, and let $P$ be a (possibly empty) path containing all vertices of $K \setminus \{k, k_1, k_2, k_3\}$. Then the concatenation of $P$ and the path $k, s_1, k_2, s_3, k_1, s_2, k_3$ is a Hamiltonian cycle of $G$. This concludes the proof. 
	\end{proof}
	
	\section{Triangle-free graphs}
	
	In this section, we prove Theorem \ref{thm:maintf}. 
	
	\begin{lemma} \label{lem:tf1}
	    Thetas, closed thetas, and triangle-free wheels are $2$-connected graphs with no Hamiltonian cycle. 
	\end{lemma}
	\begin{proof}
	    Again, $2$-connectivity can be checked easily. Thetas and closed thetas have no Hamiltonian cycles by Lemma \ref{lem:tough}, letting $X$ be the set of the ends of the (closed) theta. For a triangle-free wheel $H = (W, v)$, note that every edge $e$ of $W$ contains a vertex of degree two in $H$, and therefore every Hamiltonian cycle of $H$ contains $e$. It follows that every Hamiltonian cycle contains all edges of $W$; but these edges form a cycle that does not contain $v$, and hence no Hamiltonian cycle exists. 
	\end{proof}
	
	We assume that the reader is familar with standard definitions for graph minors and planar graphs. A \emph{model} of graph~$H$ in graph~$G$ is a collection of disjoint sets $(A_h)_{h \in V(H)}$ such that $G[A_h]$ is connected for all $h \in V(H)$, and for every edge $e = hh' \in E(H)$, there is at least one edge between $A_h$ and $A_{h'}$ in $G$. We say that graph $G$ contains $H$ \emph{as a minor} (or \emph{contains an $H$-minor}) if $G$ contains a model of $H$. A graph is \emph{outerplanar} if it has a planar embedding with all vertices incident with the outer face. 
	
	\begin{theorem}[\cite{chartrand-harary}] \label{thm:outer}
		A graph is outerplanar if and only if it has no $K_{2,3}$-minor and no $K_{4}$-minor.
	\end{theorem}
	
	\begin{lemma}[\cite{chartrand-harary}]\label{lem:outerham}
		Every $2$-connected outerplanar graph is Hamiltonian. 
	\end{lemma}

	In view of Lemma \ref{lem:tf1}, in order to prove Theorem \ref{thm:maintf}, it suffices to prove the following. 

\joe{
	\begin{theorem}
Let $G$ be a triangle-free $2$-connected graph with no induced
subgraph that is isomorphic to a theta, a closed theta, or a wheel.
Then $G$ has a Hamiltonian cycle.
	\end{theorem}
}

	\begin{proof}
		We first prove: 
		
		\sta{\label{k4} $G$ has no $K_4$-minor.}
		Suppose for a contradiction that $G$ contains a model of $K_4$ with sets $A_1, A_2, A_3, A_4$. 

\joe{
We first construct an induced cycle in $G[A_1 \cup A_2 \cup A_3]$
with at least one vertex from each of $A_1, A_2, A_3$. Let $P$ be
a shortest path in $G[A_1 \cup A_2]$ from $A_1 \cap N(A_3)$ to $A_2
\cap N(A_3)$. Let $x \in A_1$ and $y \in A_2$ be the ends of $P$.
Let $Q$ be a shortest path from $N(x) \cap A_3$ to $N(y) \cap A_3$.
Clearly, each of $P, Q$ is an induced path of $G$.  Moreover, from
the choice of $P$, no vertex of $P \setminus \{x,y\}$ has a neighbour
in $V(Q) \subseteq A_3$, and from the choice of $Q$, each of $x, y$ has
a unique neighbour in $Q$.
Thus, $C = G[V(P) \cup V(Q)]$ is an induced cycle.
}

		If there is a vertex in $G \setminus C$ that has at least two neighbours in $C$, then $G$ contains a theta or a triangle-free wheel, a contradiction. 
		
		Now, let $x \in A_4$. Let $P_i$ be a shortest path in $G[A_i \cup A_4]$ from $x$ to $V(C)$ for $i \in \{1, 2, 3\}$. For $i \in \{1,2,3\}$, let $y_i \in V(P_i) \cap V(C)$. Since $G$ is triangle-free, not all of $y_1, y_2, y_3$ are pairwise adjacent; by symmetry, say $y_1$ and $y_2$ are non-adjacent. Then $G[V(P_1) \cup V(P_2)]$ contains a path between two non-adjacent vertices of $C$; we let $R$ be a minimum-length path such that the ends $a$ and $b$ of $R$ are non-adjacent vertices of $C$, and $V(R) \setminus \{a, b\} \subseteq V(G) \setminus V(C)$. Let $a', b'$ be the neighbours of $a$ and $b$ in $R$, respectively. If no vertex in $V(R) \setminus \{a, b\}$ has a neighbour in $V(C) \setminus \{a, b\}$, then the graph $G[V(R) \cup V(C)]$ is a theta, a contradiction. 
		
\joe{
Let $w \in V(C) \setminus \{a, b\}$ be a vertex that has a neighbour
$y$ in $V(R) \setminus \{a, b\}$.  Then, $y \in V(R) \setminus
\{a,b,a',b'\}$, because $a'$, respectively, $b'$ cannot have two
(or more) neighbours in $C$.
Therefore, from the choice of $R$, it follows that $w$ is adjacent
to both $a$ and $b$. If there is a vertex $z$ in $V(C) \setminus
\{a, b, w\}$ such that $z$ has a neighbour $y'$ in $V(R) \setminus
\{a, b\}$, then again $z$ is adjacent to both $a$ and $b$, and
therefore non-adjacent to $w$. But now the path obtained from $z,
y'$, the subpath of $R$ from $y'$ to $y$, $y, w$ is shorter than
$R$, contradicting the choice of $R$.  Hence, there are no edges
between $V(R) \setminus \{a,b\}$ and $V(C) \setminus \{a,b,w\}$.
}

\joe{
It follows that $C' = G[(V(C) \cup V(R)) \setminus \{w\}]$ is an
induced cycle, and $x$ has at least three neighbours in it, namely,
$a, b, y$; thus, $G$ contains an induced triangle-free wheel, a
contradiction.  This proves \eqref{k4}.
}

		\sta{\label{k23} $G$ has no $K_{2,3}$-minor.}
		
		Suppose that $G$ has a $K_{2,3}$-minor. It follows that $G$ contains a subdivision of $K_{2,3}$ as a (not necessarily induced) subgraph. Choose such a subgraph $H$ with as few vertices as possible, and let $u$ and $v$ be the two vertices of degree three in $H$. Let the three paths between $u$ and $v$ be $P_1, P_2, P_3$. 
		
		Note that each of the paths $P_1, P_2, P_3$ is induced in $G$ except for possibly the edge $uv$, otherwise, we can replace it with a shorter path. Since $G$ has no $K_4$-minor by \eqref{k4}, it follows that there are no edges between $V(P_1) \setminus \{ u, v \}, V(P_2) \setminus \{u, v\}$ and $V(P_3) \setminus \{u, v\}$. 
		
		If $uv$ is an edge, then $G[V(H)]$ is a closed theta. If not, then $G[V(H)]$ is a theta. This is a contradiction, and proves \eqref{k23}. \medskip
		
		By \eqref{k4}, \eqref{k23}, and Theorem \ref{thm:outer}, it follows that $G$ is outerplanar. Now, the result follows from Lemma \ref{lem:outerham}.  
	\end{proof}
	
	\bibliography{haminducedsubgraphs}
\end{document}